# ASYMPTOTIC NORMALITY OF A NONPARAMETRIC ESTIMATOR OF SAMPLE COVERAGE

By Cun-Hui Zhang[1] and Zhiyi Zhang

*Rutgers University and University of North Carolina at Charlotte*

This paper establishes a necessary and sufficient condition for the asymptotic normality of the nonparametric estimator of sample coverage proposed by Good [*Biometrica* **40** (1953) 237–264]. This new necessary and sufficient condition extends the validity of the asymptotic normality beyond the previously proven cases.

**1. Introduction.** Suppose that a random sample of size $n$ is drawn (with replacement) from a population of infinitely many species. Let $X_i(n)$ be the frequency of the $i$th species in the sample. Let $\mathbf{p}_n = (p_{in}, i \geq 1)$ with $\sum_{i=1}^{\infty} p_{in} = 1$ and $P_n$ be probability measures under which the $i$th species has probability $p_{in}$ of being sampled. The infinite sequence $\mathbf{X}(n) = (X_i(n), i \geq 1)$ can be viewed as a multinomial $(n, \mathbf{p}_n)$ vector under $P_n$. For all integers $m \geq 1$

$$P_n\{X_i(n) = x_i, i = 1, \ldots, m\} = \frac{n!(1 - \sum_{i=1}^{m} p_{in})^{n-x_1-\cdots-x_m} \prod_{i=1}^{m} p_{in}^{x_i}}{(n - x_1 - \cdots - x_m)! x_1! \cdots x_m!}.$$

Let $Q_n$ be the total probability of unobserved species and $F_j(n)$ be the total number of species represented $j$ times in the sample. These random variables can be written as

$$(1.1) \quad Q_n = \sum_{i=1}^{\infty} p_{in} \delta_{i0}(n), \qquad F_j(n) = \sum_{i=1}^{\infty} \delta_{ij}(n), \qquad \delta_{ij}(n) = I\{X_i(n) = j\}.$$

Received April 2008; revised September 2008.

[1]Supported in part by NSF Grants DMS-05-04387, DMS-06-04571 and DMS-08-04626 and NSA Grant MDS 904-02-1-0063.

*AMS 2000 subject classifications.* Primary 62f10, 62F12, 62G05, 62G20; secondary 62F15.

*Key words and phrases.* Sample coverage, Turing's formula, asymptotic normality.

Good [10], while attributing an essential element of his proposal to A. M. Turing, carefully developed and studied the estimation of $Q_n$ by

$$(1.2) \qquad \widehat{Q}_n = \frac{F_1(n)}{n}.$$

The total proportion of the species not represented in the sample $Q_n$ and its estimate $\widehat{Q}_n$ have many interesting applications. For examples, Efron and Thisted [4] and Thisted and Efron [19] discuss two applications related to Shakespeare's general vocabulary and authorship of a poem; Good and Toulmin [11] and Chao [1], among many others, discuss the probability of discovering new species of animals in a population; and, more recently, Mao and Lindsay [15] study a genomic application in gene-categorization, and Zhang [20] considers applications to network species and data confidentiality problems. In addition, many authors have written about the statistical properties of $\widehat{Q}_n$. Among others, Harris [12, 13], Robbins [17], Starr [18], Holst [14], Chao [2], Esty [5, 6, 7, 8, 9] and Chao and Lee [3] are frequently referenced. However, of special relevance to the issue of concern here is Esty [6], in which the asymptotic distributional behavior of the coverage estimate under infinite dimensional probability vectors is discussed. Esty [6] gives a sufficient condition for the asymptotic normality of a $\sqrt{n}$-normalized coverage estimate. More specifically, Esty [6] proved that

$$(1.3) \qquad \lim_{n \to \infty} P_n\{Z_n \le t\} = P\{N(0,1) \le t\},$$

where

$$Z_n = \frac{n(\widehat{Q}_n - Q_n)}{\{E_n F_1(n)(1 - E_n F_1(n)/n) + 2E_n F_2(n)\}^{1/2}}$$

for all real $t$ under the sufficient condition

$$(1.4) \qquad E_n F_1(n)/n \to c_1 \in (0,1), \qquad E_n F_2(n)/n \to c_2 \ge 0.$$

Esty [6] also proved that (1.4) implies

$$(1.5) \qquad \frac{n(\widehat{Q}_n - Q_n)}{\{F_1(n)(1 - F_1(n)/n) + 2F_2(n)\}^{1/2}} \xrightarrow{\mathrm{D}} N(0,1)$$

under $P_n$.

In this paper, we extend the result of Esty [6] by establishing a necessary and sufficient condition for the asymptotic normality of the sample coverage. The family of distributions under the condition of this paper includes that of Esty [6] as a proper subset.

There are three sections in the remainder of the paper. The main results and proofs are given in Section 2. Several examples, including a few cases satisfying and a few cases not satisfying the new necessary and sufficient condition of the paper and a genomic application, are given in Section 3. The proofs of several lemmas are included in the Appendix.

**2. Main results and proofs.**

2.1. *Main results.* Define

$$
(2.1) \qquad s_{\lambda n}^2 = \sum_{i=1}^{\infty} [\lambda p_{in} e^{-\lambda p_{in}} + (\lambda p_{in})^2 e^{-\lambda p_{in}}], \qquad s_n = s_{nn}.
$$

Since $E_n F_j(n) = \sum_{i=1}^{\infty} \binom{n}{j} p_{in}^j (1 - p_{in})^{n-j}$ and $(1 - p_{in})^n \approx e^{-np_{in}}$, $s_n^2$ is an approximation of $E_n F_1(n) + 2E_n F_2(n)$.

THEOREM 1. *Let $\widehat{Q}_n = F_1(n)/n$ be the Good estimate of sample coverage $Q_n$ as in (1.2) and (1.1). Let $s_n$ be as in (2.1). Suppose that*

$$
(2.2) \qquad \limsup_{n \to \infty} E_n F_1(n)/n < 1.
$$

*Then, the central limit theorem (1.3) holds if and only if both*

$$
(2.3) \qquad E_n F_1(n) + E_n F_2(n) \to \infty
$$

*and the Lindeberg condition*

$$
(2.4) \qquad s_n^{-2} \sum_{i=1}^{\infty} (np_{in})^2 e^{-np_{in}} I\{np_{in} > \varepsilon s_n\} \to 0 \qquad \forall \varepsilon > 0
$$

*hold. In this case, (1.5) holds and*

$$
(2.5) \qquad \lim_{n \to \infty} P_n \left\{ \left| \frac{\widehat{Q}_n}{Q_n} - 1 \right| > \varepsilon \right\} = 0 \qquad \forall \varepsilon > 0.
$$

*Moreover, if (1.5) holds, then (2.3) and (2.4) imply each other.*

COROLLARY 1. *If (2.2) and (2.3) hold, then (1.3), (1.5) and (2.4) are all equivalent.*

REMARK 1. If $p_{in} = p_i$ do not depend on $n$ (under a fixed probability measure $P_n = P$), then $E_n F_1(n)/n \to 0$ always holds. In this case, Esty's [6] theorem is not applicable.

REMARK 2. We call (2.4) the Lindeberg condition, since it is equivalent to the standard Lindeberg condition when the sample size is a Poisson variable with mean $n$. Due to

$$
\sum_{i=1}^{\infty} (np_{in})^2 e^{-np_{in}} I\{np_{in} \geq M\}
$$

$$
\leq \sum_{j=0}^{\infty} M 2^{j+1} e^{-M2^j} \sum_{i=1}^{\infty} np_{in} I\{M2^j \leq np_{in} < M2^{j+1}\}
$$

$$
= O(1) nM e^{-M}
$$

with $M = \varepsilon s_n$, the Lindeberg condition (2.4) holds if $s_n/\log n \to \infty$.

REMARK 3. We prove, in Lemma 1 below, that $E_n F_1(n) + 2E_n F_2(n)$ and $s_n^2$ are within an infinitesimal fraction of each other if one of these quantities are bounded away from zero. Thus, condition (2.3) holds if and only if $s_n^2 \to \infty$.

REMARK 4. Theorem 1 is proved using Poisson approximation. The only case not covered is $E_n F_1(n)/n \to 1$, where the Poisson approximation fails and Esty's theorem does not apply.

THEOREM 2. Suppose (2.4) holds and $E_n F_1(n) \to c^* \in (0, \infty)$. Then, $E_n F_2(n) \to 0$,

$$E_n(nQ_n - c^*)^2 \to 0, \qquad n\widehat{Q}_n = F_1(n) \xrightarrow{\mathrm{D}} N_{c^*}$$

under $P_n$, where $N_{c^*}$ is a certain Poisson variable with mean $c^*$.

2.2. *Poisson approximation and proofs of theorems.* Suppose the population is sampled sequentially, so that $\mathbf{X}(m) - \mathbf{X}(m-1)$, $m \geq 1$, are i.i.d. multinomial $(1, \mathbf{p}_n)$ under $P_n$. Define

$$(2.6) \qquad \xi_n = \sum_{i=1}^{\infty} \{\delta_{i1}(n) - np_{in}\delta_{i0}(n)\} = n(\widehat{Q}_n - Q_n).$$

Let $N_\lambda$ be a Poisson process independent of $\{\mathbf{X}(m), m \geq 1\}$ with $E_n N_\lambda = \lambda$. Define

$$(2.7) \qquad \zeta_{\lambda n} = \sum_{i=1}^{\infty} Y_{i\lambda n}, \qquad Y_{i\lambda n} = \delta_{i1}(N_\lambda) - \lambda p_{in}\delta_{i0}(N_\lambda).$$

Under probability $P_n$, $\{X_i(N_\lambda), i \geq 1\}$ are independent Poisson variables with means $\lambda p_{in}$, so that $\{Y_{i\lambda n}, i \geq 1\}$ are independent zero-mean variables with

$$(2.8) \qquad \begin{aligned} E_n Y_{i\lambda n}^2 &= \sigma_{i\lambda n}^2 = \lambda p_{in} e^{-\lambda p_{in}} + (\lambda p_{in})^2 e^{-\lambda p_{in}}, \\ E_n \zeta_{\lambda n}^2 &= \sum_{i=1}^{\infty} \sigma_{i\lambda n}^2 = s_{\lambda n}^2. \end{aligned}$$

THEOREM 3. Suppose $\lambda = \lambda_n \to \infty$. Then,

$$(2.9) \qquad \zeta_{\lambda n}/s_{\lambda n} \xrightarrow{\mathrm{D}} N(0,1),$$

if and only if both $s_{\lambda n} \to \infty$ and

$$(2.10) \qquad s_{\lambda n}^{-2} \sum_{i=1}^{\infty} (\lambda p_{in})^2 e^{-\lambda p_{in}} I\{\lambda p_{in} > \varepsilon s_{\lambda n}\} \to 0 \qquad \forall \varepsilon > 0.$$

PROOF OF THEOREM 3. By the Lindeberg–Feller central limit theorem, (2.9) holds if and only if

$$\max_{i \geq 1} \sigma_{i\lambda n}^2 / s_{\lambda n}^2 = \max_{i \geq 1} s_{\lambda n}^{-2} [\lambda p_{in} e^{-\lambda p_{in}} + (\lambda p_{in})^2 e^{-\lambda p_{in}}] \to 0, \tag{2.11}$$

and the standard Lindeberg condition holds in the form

$$s_{\lambda n}^{-2} \sum_{i=1}^{\infty} E Y_{i\lambda n}^2 I\{|Y_{i\lambda n}| > \varepsilon s_{\lambda n}\} \to 0 \qquad \forall \varepsilon > 0. \tag{2.12}$$

Since $\delta_{ij}(N_\lambda)$ are 0–1 variables and $Y_{i\lambda n}^2 = \delta_{i1}(N_\lambda) + (\lambda p_{in})^2 \delta_{i0}(N_\lambda)$,

$$2^{-1} Y_{i\lambda n}^2 I\{Y_{i\lambda n}^2 > 2(\varepsilon s_{\lambda n})^2\}$$

$$\leq \delta_{i1}(N_\lambda) I\{1 > \varepsilon s_{\lambda n}\} + (\lambda p_{in})^2 \delta_{i0}(N_\lambda) I\{\lambda p_{in} > \varepsilon s_{\lambda n}\},$$

which is no greater than $Y_{i\lambda n}^2 I\{|Y_{i\lambda n}| > \varepsilon s_{\lambda n}\}$. Thus, (2.12) is equivalent to

$$s_{\lambda n}^{-2} \sum_{i=1}^{\infty} [\lambda p_{in} e^{-\lambda p_{in}} I\{1 > \varepsilon s_{\lambda n}\} + (\lambda p_{in})^2 e^{-\lambda p_{in}} I\{\lambda p_{in} > \varepsilon s_{\lambda n}\}] \to 0 \tag{2.13}$$

$$\forall \varepsilon > 0.$$

If $s_{\lambda n} \to \infty$, then (2.10) implies (2.13) immediately and (2.11) via $(\lambda p_{in})^j e^{-\lambda p_{in}} \leq j!, j = 1, 2$.

It remains to prove that (2.11) and (2.13) together imply $s_{\lambda n} \to \infty$ and (2.10). In fact, (2.11) is not even needed. If $s_{\lambda n} \leq M$ along a subsequence, then, for $\varepsilon < 1/M$,

$$s_{\lambda n}^2 \leq \sum_{i=1}^{\infty} [2\lambda p_{in} e^{-\lambda p_{in}} + (\lambda p_{in})^2 e^{-\lambda p_{in}} I\{\lambda p_{in} > 1 > \varepsilon s_{\lambda n}\}]$$

$$\leq 2 \sum_{i=1}^{\infty} [\lambda p_{in} e^{-\lambda p_{in}} I\{1 > \varepsilon s_{\lambda n}\} + (\lambda p_{in})^2 e^{-\lambda p_{in}} I\{\lambda p_{in} > \varepsilon s_{\lambda n}\}],$$

so that (2.13) fails. Thus, (2.13) implies $s_{\lambda n} \to \infty$. This completes the proof, since (2.13) implies (2.10) immediately. □

We prove Theorems 1 and 2 via Theorem 3 and the Poisson approximation

$$\frac{\xi_n - \zeta_{nn}}{s_n} = o_{P_n}(1). \tag{2.14}$$

We need three lemmas.

LEMMA 1. (i) *Let* $s_n^2$ *be as in* (2.1). *For* $\varepsilon/n \leq 1/4$,

$$(1 - 1/n) e^{-\varepsilon} s_n^2 - n^2 e^{-\sqrt{\varepsilon n}} \leq E_n F_1(n) + 2 E_n F_2(n) \leq e^{2\varepsilon} s_n^2 + n(n+1) e^{-(n-2)\varepsilon}.$$

*Consequently, if* $\liminf_n \min\{s_n^2, E_n F_1(n) + E_n F_2(n)\} > 0$, *then*

$$\{E_n F_1(n) + 2 E_n F_2(n)\} / s_n^2 \to 1.$$

(ii) *Let* $s_{\lambda n}^2$ *and* $s_n^2$ *be as in* (2.1). *For all* $\lambda' < \lambda$ *and* $\varepsilon > 0$,

$$(2.15) \quad (\lambda'/\lambda)^2 s_{\lambda n}^2 \le s_{\lambda' n}^2 \le e^\varepsilon s_{\lambda n}^2 + \lambda(1+\lambda)\exp(-\lambda'\varepsilon/(\lambda-\lambda')).$$

*Consequently,* $s_{\lambda_n n}^2 = (1+o(1))s_n^2$ *if* $n^2 e^{-\varepsilon n/|\lambda_n - n|} = o(s_n^2)$ *for all* $\varepsilon > 0$ *and* $\lambda_n/n \to 1$.

LEMMA 2. *Let* $\zeta_{\lambda n}$ *be as in* (2.7). *Then,*

$$E_n \max_{\lambda \le t \le \lambda+\Delta} |\zeta_{tn} - \zeta_{\lambda n}|$$

$$\le 2\left\{\sum_{i=1}^\infty \lambda p_{in}(1+\lambda p_{in})e^{-\lambda p_{in}}(1-e^{-\Delta p_{in}})\right\}^{1/2} + 2\sum_{i=1}^\infty \Delta p_{in}e^{-\lambda p_{in}}.$$

LEMMA 3. *If* $\liminf_n s_n^2 > 0$ *and* $s_n^2/n = o(1)$, *then* (2.14) *holds.*

PROOF OF THEOREM 2. It follows, from (1.1) and (2.4), that $2E_n F_2(n)$ is bounded by

$$\sum_{i=1}^\infty 2\binom{n}{2} p_{in}^2 (1-p_{in})^{n-2} \le \sum_{i=1}^\infty (np_{in})^2\{(1-p_{in})^{n-1}+p_{in}\}I\{np_{in} \le \varepsilon s_n\}$$

$$(2.16) \qquad\qquad\qquad + \sum_{i=1}^\infty (np_{in})^2 e^{-(n-2)p_{in}} I\{np_{in} > \varepsilon s_n\}$$

$$\le \varepsilon s_n E_n F_1(n) + (\varepsilon s_n)^2 + o(s_n^2),$$

so that, due to $E_n F_1(n) = O(1)$, $s_n^2 = O(1)$ by Lemma 1(i). Thus, by (2.1) and (2.4),

$$(2.17) \quad \sum_{i=1}^\infty (np_{in})^2 e^{-np_{in}} \le \sum_{i=1}^\infty (np_{in})^2 e^{-np_{in}} I\{np_{in} > \varepsilon s_n\} + \varepsilon s_n^3 \to 0$$

as $n \to \infty$ and then $\varepsilon \to 0+$. Since $E_n \delta_{ij}(n) = \binom{n}{j} p_{in}^j (1-p_{in})^{n-j}$, (2.17) implies

$$0 \le E_n\{F_1(n) - nQ_n\} = \sum_{i=1}^\infty np_{in}\{(1-p_{in})^{n-1} - (1-p_{in})^n\}$$

$$\le e\sum_{i=1}^\infty np_{in}^2 e^{-np_{in}} \to 0,$$

so that $nE_n Q_n \to c^*$. Since $\{\delta_{i0}(n), i \ge 1\}$ have negative correlation, (2.17) also implies

$$\text{Var}_n(nQ_n) \le \sum_{i=1}^\infty \text{Var}(np_{in}\delta_{i0}(n)) \le \sum_{i=1}^\infty (np_{in})^2 e^{-np_{in}} \to 0.$$

Thus, $E_n(nQ_n - c^*)^2 \to 0$. Similarly, $E_n F_2(n) \le (e^2/2) \sum_{i=1}^\infty (np_{in})^2 e^{-np_{in}} \to 0$.

Let $\widetilde{Q}_n = \sum_{i=1}^\infty p_{in} \delta_{i0}(N_n)$. By (2.17), $\mathrm{Var}_n(n\widetilde{Q}_n) = \sum_{i=1}^\infty (np_{in})^2 e^{-np_{in}} = o(1)$. By (2.17) and then Lemma 1(i), $nE\widetilde{Q}_n = \sum_{i=1}^\infty np_{in} e^{-np_{in}} = s_n^2 + o(1) = c^* + o(1)$. These imply $n\widetilde{Q}_n = c^* + o_{P_n}(1)$. Thus, by Lemma 3,

$$F_1(n) - F_1(N_n) = \xi_n + nQ_n - (\zeta_{nn} + n\widetilde{Q}_n) = \xi_n - \zeta_{nn} + o_{P_n}(1) = o_{P_n}(1).$$

Since $F_1(N_n) = \sum_{i=1}^\infty \delta_{i1}(N_n)$ are independent Bernoulli variables with uniformly small probabilities $E_n \delta_{i1}(N_n) = np_{in} e^{-np_{in}} \le \{\sum_{i=1}^n (np_{in})^2 e^{-np_{in}}\}^{1/2} = o(1)$, $F_1(n) = F_1(N_n) + o_{P_n}(1)$ converges in distribution to a Poisson variable with mean $E_n F_1(N_n) = nE\widetilde{Q}_n \to c^*$. $\quad\square$

PROOF OF THEOREM 1. Assume, without loss of generality, that

$$E_n F_j(n)/n \to c_j, \qquad j = 1, 2, \qquad E_n F_1(n) + 2E_n F_2(n) \to c^*,$$

with $c_1 \in [0, 1)$, $c_2 \in [0, 1]$ and $c^* \in [0, \infty]$ (taking subsequence if necessary).

Case 1. $c_1 > 0$. It follows from the theorem of Esty [6] that (1.3) holds. Moreover, since $s_n^2/n \to c_1 + 2c_2 > 0$ by Lemma 1(i), (2.4) holds as in Remark 2. Thus, (1.3), (2.3) and (2.4) all hold.

Case 2. $c_1 = c^* = 0$. Since $E_n F_1(n) \to 0$ and $Z_n \le 0$ for $F_1(n) = 0$,

$$P_n(Z_n \le 0) \ge P_n(F_1(n) = 0) \to 1.$$

Thus, (1.3) does not hold. Similarly, (1.5) does not hold. Since $c^* = 0$, (2.3) does not hold.

Case 3. $c_1 = 0 < c^*$. By (1.1), $2E_n F_2(n)/(n-1)$ is bounded by

$$\sum_{i=1}^\infty np_{in}^2(1-p_{in})^{n-2} \le \frac{M}{1-M/n} \sum_{i=1}^\infty p_{in}(1-p_{in})^{n-1} + \sup_{p \ge M/n} np(1-p)^{n-2}.$$

Since $\sum_{i=1}^\infty p_{in}(1-p_{in})^{n-1} = E_n F_1(n)/n \to c_1 = 0$, we find $E_n F_2(n)/n \to 0 = c_2$, which then implies $s_n^2/n \to 0$ by Lemma 1(i). In addition, Lemma 1(i) implies $\{E_n F_1(n) + 2E_n F_2(n)\}/s_n^2 \to 1$, so that $s_n^2 \to c^* > 0$. Thus, (2.14) holds by Lemma 3, and (1.3) holds if and only if $\zeta_{nn}/s_n \to N(0,1)$ in view of (2.6). Therefore, by Theorem 3 with $\lambda = n$, (1.3) holds if and only if both (2.3) and (2.4) hold.

We have proved the first assertion of the theorem, since (1.3) holds if and only if both (2.3) and (2.4) hold in all the three cases. It remains to prove that (1.3) implies (1.5) and (2.5), and that (2.3) and (2.4) are equivalent under (1.5).

We first prove the equivalence of (1.3) and (1.5) under (2.3). For fixed $(j, n)$, $\delta_{ij}(n)$ are Bernoulli variables with $\mathrm{Cov}_n(\delta_{ij}(n), \delta_{i'j}(n)) \le 0$, so that $\mathrm{Var}_n(F_j(n)) \le E_n F_j(n)$ and

$$\mathrm{Var}_n(F_1(n) + 2F_2(n)) \le 2\{E_n F_1(n) + 4E_n F_2(n)\}.$$

Since $E_n F_1(n) + 2E_n F_2(n) \to \infty$, $\{F_1(n) + 2F_2(n)\}/\{E_n F_1(n) + 2E_n F_2(n)\} \to$ 1 in $P_n$ by the above inequality. Similarly, $F_1^2(n)/n = (1 + o_{P_n}(1))\{E_n F_1(n)\}^2/n$. Moreover, since $\{E_n F_1(n)\}^2/n = (c_1 + o(1)) E_n F_1(n)$ with $c_1 < 1$, $E_n F_1(n)\{1 - E_n F_1(n)/n\} + 2F_2(n)$ is of the same order as $E_n F_1(n) + 2E_n F_2(n)$. Thus, (1.3) and (1.5) are equivalent under (2.3).

Assume (1.3) holds. Since (2.3) holds, (1.5) holds. Since the Lindeberg (2.4) holds,

$$(2.18) \qquad 2E_n F_2(n) = o(s_n) E_n F_1(n) + o(s_n^2)$$

by (2.16). Thus, (2.3) and Lemma 1(i) provide

$$s_n^2 = (1 + o(1))\{E_n F_1(n) + 2E_n F_2(n)\} = (1 + o(s_n)) E_n F_1(n) + o(s_n^2) \to \infty,$$

which implies $s_n = o(1) E_n F_1(n)$ and $\mathrm{Var}_n(F_1(n)) \le E_n F_1(n) \to \infty$. Consequently, $s_n = o_{P_n}(F_1(n))$, and then, by (1.3), $nQ_n - n\widehat{Q}_n = O_{P_n}(s_n) = o_{P_n}(F_1(n)) = o_{P_n}(n\widehat{Q}_n)$. Thus, (1.3) implies (2.5) as well as (1.5).

Now, we assume (1.5). If (2.3) holds, then (1.3) holds due to its equivalence to (1.5), so that (2.4) must hold. It remains to prove (2.3); that is, $c^* = \infty$ under (2.4). Since (1.5) holds, Case 2 is ruled out, so $c^* > 0$. If $0 < c^* < \infty$, Lemma 1(i) implies $s_n^2 = (1 + o(1))\{E_n F_1(n) + 2E_n F_2(n)\} = O(1)$, and then (2.18) implies $E_n F_2(n) = o(1)$, so that $E_n F_1(n) \to c^*$. Thus, by Theorem 2, $0 < c^* < \infty$ would imply the convergence of $\sqrt{c^*} Z_n$ in distribution to $N_{c^*} - c^*$ and the convergence of $F_1(n)(1 - F_1(n)/n) + 2F_2(n)$ to $N_{c^*}$. This is impossible since (1.5) holds. Hence, $c^* = \infty$. □

## 3. Examples.

We provide three theoretical examples and describe one real application. In all theoretical examples, we define $p_{in} \propto p_n(i)$ with $\int_0^\infty p_n(x)\,dx = 1$. The density functions $p_n(x)$ are decreasing in $x > 0$ and sufficiently regular to allow the following approximations within an infinitesimal fraction:

$$(3.1) \qquad \begin{aligned} E_n F_1(n) &\approx \int_0^\infty n p_n(x) e^{-n p_n(x)}\,dx, \\ s_n^2 &\approx \int_0^\infty n p_n(x)\{1 + n p_n(x)\} e^{-n p_n(x)}\,dx. \end{aligned}$$

EXAMPLE 1 (*Fixed discrete Paretos*). In this example, Theorem 1 provides the asymptotic normality, but the Esty's [6] condition $E_n F_1(n)/n \to c_1 \in (0, 1)$ does not hold. Let $p_n(x) = p(x) = a/(x + 1)^b$ with $a > 0$ and $b > 1$. Condition (2.2) is satisfied, since $E_n F_1(n)/n \approx \int_0^\infty p(x) e^{-n p(x)}\,dx \to 0$. For large $n$, changing variable $t = np(x) = na/(x + 1)^b$ yields

$$E_n F_1(n) \approx -\int_0^{na} t e^{-t}\,d(na/t)^{1/b} \approx \frac{(na)^{1/b}}{b} \int_0^\infty t^{-1/b} e^{-t}\,dt \propto n^{1/b},$$

so that (2.3) holds and $s_n/\log n \to \infty$ by Lemma 1(i). It follows that (2.4) holds by Remark 2. Thus, the central limit theorems (1.3) and (1.5) both hold by Theorem 1.

EXAMPLE 2 (*Dynamic discrete exponentials*). In this example, (2.3) and (2.4) are equivalent. Let $p_n(x) = a_n^{-1}e^{-x/a_n}$ with $a_n/n \le M < \infty$. Let $t = np_n(x)$. By (3.1),

$$\frac{E_nF_1(n)}{n} \approx n^{-1}\int_0^{n/a_n} te^{-t}\,d(a_n\log t) = \int_0^1 e^{-yn/a_n}\,dy < 1,$$

so that (2.2) holds. Similarly, $s_n^2 \approx a_n\int_0^{n/a_n}\{1+t\}e^{-t}\,dt$ by (3.1), so that $s_n^2$ is of the order $a_n$. Moreover, the Lindeberg condition (2.4) is equivalent to

$$o(1) = \frac{1}{a_n}\int_{np_n(x)>\varepsilon\sqrt{a_n}}\{np_n(x)\}^2 e^{-np_n(x)}\,dx = \int_{\varepsilon\sqrt{a_n}<t<n/a_n} te^{-t}\,dt,$$

which holds if and only if $s_n^2 \sim a_n \to \infty$, if and only if (2.3) holds by Lemma 1(i).

EXAMPLE 3 (*Dynamic two-step functions*). This example demonstrates that the three conditions of Theorem 1 are not redundant. Let $a_{jn} \to \infty$ and $w_{1n} + w_{2n} = 1$ with $w_{1n}/a_{1n} \ge w_{2n}/a_{2n} \ge 0$. Set $p_n(x) = \sum_{j=1}^2 w_{jn}a_{jn}^{-1}I\{0 < (-1)^j(x-a_{1n}) \le a_{jn}\}$. By (3.1),

$$E_nF_1(n) \approx n\sum_{j=1}^2 w_{jn}e^{-b_{jn}},$$

$$s_n^2 \approx n\sum_{j=1}^2 w_{jn}(1+b_{jn})e^{-b_{jn}}, \qquad b_{jn} = nw_{jn}/a_{jn}.$$

Moreover, the Lindeberg condition (2.4) holds if and only if

$$\frac{n}{s_n^2}\sum_{j=1}^2 w_{jn}b_{jn}e^{-b_{jn}}I\{b_{jn} > \varepsilon s_n\} \to 0 \qquad \forall \varepsilon > 0.$$

Case 1. $w_{1n} = 1$ and $b_{1n} \not\to 0$. The $p_n(x)$ are uniform densities in $(0, a_{1n})$. Condition (2.2) holds, since $E_nF_1(n)/n \approx e^{-b_{1n}} \not\to 1$. Since $1 + b_{1n}$ is of the same order as $b_{1n}$, (2.4) holds if and only if $b_{1n}/s_n \to 0$, so that (2.4) implies (2.3). Let $b_{1n} = \log n - \log\log n$. We find $s_n^2 \approx (1+b_{1n})\log n \approx b_{1n}^2 \to \infty$. Thus, both (2.2) and (2.3) hold but (2.4) does not.

Case 2. $w_{1n} = 1$ and $b_{1n} \to 0$. The $p_n(x)$ are still uniform. Since $E_nF_1(n)/n \approx e^{-b_{1n}} \to 1$, (2.2) does not hold. On the other hand, $s_n^2 \approx n(1+b_{1n})e^{-b_{1n}} \to \infty$ and $b_{1n}/s_n \to 0$. Thus, both (2.3) and (2.4) hold but (2.2) does not.

Case 3. $w_{1n} = (1 - 1/n)$, $b_{1n} = 2 \log n$ and $b_{2n} \to 0$. Since $E_n F_1(n)/n = o(1)$ and $s_n^2 = o(1) + n w_{2n}(1 + o(1)) \to 1$, both (2.2) and (2.4) hold but (2.3) does not.

EXAMPLE 4 (*A genomic application*). Mao and Lindsay [15] studied a gene expression problem based on a sample of $n = 2568$ expressed sequence tags from a tomato flower cDNA library. The data came from the Institute for Genomic Research. Detailed description of the data set may also be found in Quackenbush et al. [16]. In this context, $Q_n$ is the probability that the next randomly selected expressed sequence tag will stand for a new gene. A quantification of $Q_n$ will then be an informative indicator pertaining to the depth of the sample collected thus far regarding the levels of expression of the genes in the library. For this particular data set, $n = 2568$, $F_1(n) = 1434$, $F_2(n) = 253$, $F_3(n) = 71$, $F_4(n) = 33$, $F_5(n) = 11$, $F_6(n) = 6$, $F_7(n) = 2$, $F_8(n) = 3$, $F_9(n) = 1$, $F_{10}(n) = F_{11}(n) = 1$ and $F_{12}(n) = F_{13}(n) = F_{14}(n) = F_{16}(n) = F_{23}(n) = F_{27}(n) = 1$, resulting in $\widehat{Q}_n = 0.5584$. By (1.5), the 95% confidence interval for $Q_n$ is $(0.5391, 0.5777)$, which incidentally is narrower than the 95% confidence interval produced by Mao and Lindsay [15], $(0.529, 0.580)$. Our confidence interval is not new, since it was based on an identical expression given by Esty [6]. However, we take a bit more comfort in such applications, in knowing that the validity of the confidence interval is supported by a larger family of distributions as a result of Theorem 1.

REMARK 5. The procedure introduced by Mao and Lindsay [15] is applicable to not only the total probability associated with nonrepresented genes but also that associated with genes represented with frequencies lower than a threshold. They took a different perspective to the problem from that of Esty [6] and, hence, ours. Specifically, their derivation started by directly assuming $(X_i(n), i \geq 1)$, being independent Poisson random variables with means $(\lambda_i, i \geq 1)$ which is itself an i.i.d. sample from a latent distribution. Their results are based on an asymptotical argument with the number of species (genes) approaching infinity.

## APPENDIX: PROOFS OF LEMMAS

PROOF OF LEMMA 1. (i) Since $1 - p \leq e^{-p}$,

$$E_n F_1(n) + 2 E_n F_2(n) = \sum_{i=1}^{\infty} \{ n p_{in}(1 - p_{in})^{n-1} + n(n-1) p_{in}^2 (1 - p_{in})^{n-2} \}$$

$$\leq \sum_{i=1}^{\infty} n p_{in}(1 + n p_{in}) e^{-(n-2) p_{in}}$$

$$\leq e^{2\varepsilon}s_n^2 + \sum_{i=1}^{\infty} np_{in}(1+n)e^{-(n-2)\varepsilon}.$$

Since $1 - p \geq e^{-p-p^2}$ for $0 \leq p \leq 1/2$ and $1 - p + (n-1)p \geq (1-1/n)(1-p)^2(1+np)$,

$$E_n F_1(n) + 2E_n F_2(n) = \sum_{i=1}^{\infty} np_{in}(1-p_{ni})^{n-2}(1-p_{ni}+(n-1)p_{ni})$$

$$\geq (1-1/n)\sum_{i=1}^{\infty} np_{in}(1+np_{in})e^{-np_{in}-\varepsilon}I\{np_{in}^2 \leq \varepsilon\}$$

$$\geq (1-1/n)e^{-\varepsilon}s_n^2 - n^2 e^{-\sqrt{\varepsilon n}}.$$

(ii) For all $\lambda' < \lambda$ and $\varepsilon > 0$,

$$(\lambda'/\lambda)^2 s_{\lambda n}^2 \leq s_{\lambda' n}^2$$

$$\leq \sum_{i=1}^{\infty} \lambda p_{in}(1+\lambda p_{in})e^{-\lambda' p_{in}}$$

$$\leq e^{\varepsilon}s_{\lambda n}^2 + \sum_{i=1}^{\infty} \lambda p_{in}(1+\lambda p_{in})e^{-\lambda' p_{in}}I\{(\lambda-\lambda')p_{in} > \varepsilon\}$$

$$\leq e^{\varepsilon}s_{\lambda n}^2 + \lambda(1+\lambda)\exp\big(-\lambda'\varepsilon/(\lambda-\lambda')\big).$$

This gives (2.15), and the rest follows easily. $\square$

PROOF OF LEMMA 2. Let $Y_{i\lambda n} = \delta_{i1}(N_\lambda) - \lambda p_{in}\delta_{i0}(N_\lambda)$ be as in (2.7). For $t > \lambda$,

$$\begin{aligned}
Y_{itn} - Y_{i\lambda n} &= \delta_{i1}(N_t) - tp_{in}\delta_{i0}(N_t) - \delta_{i1}(N_\lambda) + \lambda p_{in}\delta_{i0}(N_\lambda) \\
&= \delta_{i1}(N_\lambda)\{\delta_{i1}(N_t) - 1\} \\
&\quad + \delta_{i0}(N_\lambda)\{\delta_{i1}(N_t) - tp_{in}\delta_{i0}(N_t) + \lambda p_{in}\delta_{i0}(N_\lambda)\} \\
&= -Y_{i\lambda n}I\{X_i(N_t) > X_i(N_\lambda)\} \\
&\quad + \delta_{i0}(N_\lambda)\{\delta_{i1}(N_t) - (t-\lambda)p_{in}\delta_{i0}(N_t)\}.
\end{aligned}$$ (A.1)

The above identity can be verified by checking both the cases of $\delta_{i0}(N_\lambda) \in \{0,1\}$ and by noticing that $\delta_{ij}(N_\lambda)\{1-\delta_{ij}(N_t)\} = \delta_{ij}(N_\lambda)I\{X_i(N_t) > X_i(N_\lambda)\}$.

Let $T_i = \min\{t : X_i(N_t) > X_i(N_\lambda)\}$. Since $\{Y_{i\lambda n}, i \geq 1\}$ are independent variables with mean zero and independent of $\{\mathbf{X}(N_t) - \mathbf{X}(N_\lambda), t \geq \lambda\}$, by Doob's inequality for martingales,

$$E_n \max_{\lambda < t \leq \lambda+\Delta}\left[\sum_{i=1}^{\infty} Y_{i\lambda n}I\{X_i(N_t) > X_i(N_\lambda)\}\right]^2$$

$$
\begin{aligned}
&= E_n \max_{\lambda < t \le \lambda + \Delta} \left[ \sum_{T_i \le t} Y_{i\lambda n} \right]^2 \\
&\le 4 \sum_{i=1}^{\infty} E_n Y_{in}^2(\lambda) I\{X_i(N_{\lambda+\Delta}) > X_i(N_\lambda)\} \\
&= 4 \sum_{i=1}^{\infty} \lambda p_{in}(1 + \lambda p_{in}) e^{-\lambda p_{in}}(1 - e^{-\Delta p_{in}}).
\end{aligned}
$$

(A.2)

For the second term on the right-hand side of (A.1), we have

$$
\begin{aligned}
E_n &\sup_{\lambda < t \le \lambda + \Delta} \left| \sum_{i=1}^{\infty} \delta_{i0}(N_\lambda)\{\delta_{i1}(N_t) - (t - \lambda) p_{in} \delta_{i0}(N_t)\} \right| \\
&\le \sum_{i=1}^{\infty} E_n \delta_{i0}(N_\lambda)(P_n\{X_i(N_{\lambda+\Delta}) > X_i(N_\lambda)\} + \Delta p_{in}) \\
&\le \sum_{i=1}^{\infty} e^{-\lambda p_{in}} 2\Delta p_{in}.
\end{aligned}
$$

This and (A.2) yield the conclusion in view of (A.1).  □

PROOF OF LEMMA 3.   Let $t_n$ be the arrival time of the $n$th event in the Poisson process $N_\lambda$, with $N_{t_n} = n$. Since $\xi_n - \zeta_{t_n n} = (t_n - n) \sum_{i=1}^{\infty} p_{in} \delta_{i0}(n)$, we have

$$
\begin{aligned}
P_n&\{|\xi_n - \zeta_{nn}| > \varepsilon s_n\} \\
&\le P_n\{|t_n - n| > \Delta/2\} \\
&\quad + P_n\left\{ \max_{n-\Delta/2 < t < n+\Delta/2} |\zeta_n - \zeta_{tn}| + (\Delta/2) \sum_{i=1}^{\infty} p_{in} \delta_{i0}(n) > \varepsilon s_n \right\}.
\end{aligned}
$$

(A.3)

Set $\lambda = n - \Delta/2$. Since $E_n \delta_{i0}(n) = (1 - p_{in})^n \le e^{-np_{in}} \le e^{-\lambda p_{in}}$, by Lemma 2,

$$
\begin{aligned}
E_n&\left\{ \max_{n-\Delta/2 < t < n+\Delta/2} |\zeta_n - \zeta_{tn}| + (\Delta/2) \sum_{i=1}^{\infty} p_{in} \delta_{i0}(n) \right\} \\
&\le 4\left\{ \sum_{i=1}^{\infty} \lambda p_{in}(1 + \lambda p_{in}) e^{-\lambda p_{in}}(1 - e^{-\Delta p_{in}}) \right\}^{1/2} \\
&\quad + (4 + 1/2) \sum_{i=1}^{\infty} \Delta p_{in} e^{-\lambda p_{in}}.
\end{aligned}
$$

(A.4)

Since $t_n$ has the gamma$(n, 1)$ distribution, $E_n(t_n - n)^2 = n$. Thus, by (A.3)

and (A.4), (2.14) holds via the Markov inequality, provided that

$$s_n^{-2} \sum_{i=1}^{\infty} \lambda p_{in}(1+\lambda p_{in})e^{-\lambda p_{in}}(1-e^{-\Delta p_{in}}) \to 0,$$

(A.5)

$$\frac{\Delta}{s_n} \sum_{i=1}^{\infty} p_{in}e^{-\lambda p_{in}} \to 0,$$

with $n - \lambda = \Delta = M\sqrt{n} = O(\sqrt{\lambda})$ for all $0 < M < \infty$.

It remains to prove (A.5). Since $\liminf_n s_n^2 > 0$, $s_{\lambda n}/s_n \to 1$ by Lemma 1(ii). Since $s_n^2/n = o(1)$, the second part of (A.5) holds due to $(\Delta/s_n)\sum_{i=1}^{\infty} p_{in}e^{-\lambda p_{in}} \le s_{\lambda n}^2 \Delta/\lambda s_n = O(1)s_n/\sqrt{n} = o(1)$. For the first part of (A.5),

$$\sum_{i=1}^{\infty} \lambda p_{in}(1+\lambda p_{in})e^{-\lambda p_{in}}(1-e^{-\Delta p_{in}})$$

$$\le \varepsilon s_{\lambda n}^2 + \sum_{i=1}^{\infty} \lambda p_{in}(1+\lambda p_{in})e^{-\lambda p_{in}}I\{\Delta p_{in} > \varepsilon\}$$

$$\le \varepsilon s_{\lambda n}^2 + \lambda(1+\lambda)e^{-\lambda\varepsilon/\Delta} \le (1+o(1))\varepsilon s_n^2 + o(1).$$

Thus, since $\liminf_n s_n^2 > 0$, the proof is complete. $\square$

DEPARTMENT OF STATISTICS
RUTGERS UNIVERSITY
NEW BRUNSWICK, NEW JERSEY 08903
USA
E-MAIL: czhang@stat.rutgers.edu

DEPARTMENT OF MATHEMATICS
  AND STATISTICS
UNIVERSITY OF NORTH CAROLINA
  AT CHARLOTTE
CHARLOTTE, NORTH CAROLINA 28223
USA
E-MAIL: zzhang@uncc.edu